\documentclass[11pt,intlimits, draft,oneside]{amsart}

\usepackage[latin2]{inputenc}
\usepackage{latexsym}
\usepackage{amssymb}
\usepackage{amsthm}

\usepackage{color} 

\usepackage{amsmath}

\usepackage{enumerate}

\def\Eins{\mathbf{1}}
\input{cyracc.def}
\newfont{\cyrfnt}{wncyi10 at 11pt}

\newtheorem{thm}{Theorem}[section]
\newtheorem{prop}[thm]{Proposition}
\newtheorem{lemma}[thm]{Lemma}

\theoremstyle{definition}
\newtheorem{defi}[thm]{Definition}
\newtheorem{remark}[thm]{Remark}
\newtheorem{example}[thm]{Example}

\newcommand{\R}{\mathbb{R}}
\newcommand{\Q}{\mathbb{Q}}
\newcommand{\C}{\mathbb{C}}
\newcommand{\CC}{\mathcal{C}}
\newcommand{\I}{\mathcal{I}}
\newcommand{\N}{\mathbb{N}}

\newcommand{\U}{\mathcal{U}}

\newcommand{\ra}{\rangle}
\newcommand{\la}{\langle}
\newcommand{\veps}{\varepsilon}

\newcommand{\vphi}{\varphi}
\newcommand{\ol}{\overline}

\newcommand{\cont}{\text{cont}}
\newcommand{\lin}{\text{lin}}
\renewcommand{\Re}{\mathrm{Re}\,}

\def\LLL{\mathcal{L}}

\begin{document}

\title{Rigidity of contractions on Hilbert spaces}

\author{Tanja Eisner}
\address{Tanja Eisner \newline Mathematisches Institut, Universit\"{a}t T\"{u}bingen\newline Auf der Morgenstelle 10, D-72076, T\"{u}bingen, Germany}
\email{talo@fa.uni-tuebingen.de}

\renewcommand{\thefootnote}{}


\begin{abstract}
We study the asymptotic behaviour of contractive operators and
strongly continuous semigroups on separable Hilbert spaces 
using the notion of rigidity.
In particular, we show that a ``typical'' contraction $T$ 
contains the unit circle times the identity operator in the strong limit set of its powers,
while $T^{n_j}$ converges weakly to zero along a sequence $\{n_j\}$
with density one. The continuous analogue is presented for isometric and unitary $C_0$-(semi)groups. 
\end{abstract}

\maketitle



\section{Introduction}

``Good'' behaviour of the powers $T^n$ of a linear operator $T$ on a
Banach space $X$ has been studied intensively leading to many important applications. 
Here, ``good'' behaviour may mean ``stability'' in the sense that $\lim_{n\to\infty}
T^n = 0$ with respect to one of the standard operator topologies. We
refer to M\"uller \cite{mueller:2007}, Chill, Tomilov \cite{chill/tomilov:2007} and \cite{eisner-book} for a survey on these
properties. Further ``good'' properties might be  
convergence to a projection, to a periodic or compact group. 
On the other hand, it is also well-known that linear operators and their powers may
behave quite differently and, e.g., have a dense orbit $\{T^nx,n\in\N\}$ for some $x$ in
$X$.
See the recent monograph of Bayart, Matheron \cite{bayart/matheron:2009} on such hypercyclic or chaotic operators.

In this paper we look at a still different type of asymptotic behaviour. 
In the context of contractions
on separable Hilbert spaces we show that it can often happen that for
some (necessarily very large) subsequence $\{n_j\}_{j=1}^\infty\subset\N$ 
weak-$\lim_{j\to\infty} T^{n_j} = 0$ while for every $\lambda\in
\Gamma$, $\Gamma$ the unit circle, strong-$\lim_{j\to\infty} T^{m_j}
= \lambda I$
for some other sequence $\{m_j\}$ 
(this second property is called $\Gamma$-rigidity and for $\lambda=1$ rigidity). Analogous rigidity properties have been studied for measures on the
unit circle (see Nadkarni \cite[Chapter 7]{nadkarni:1998}) and for
measure theoretical dynamical systems (see e.g. Katok
\cite{katok:1985}, Nadkarni \cite[Chapter 8]{nadkarni:1998}, Goodson,
Kwiatkowski, Lemanczyk, Liardet \cite{lemanczyk-etal:1992} and
Ferenczi, Holton, Zamboni \cite{ferenczi/holton/zamboni:2005}). 
In our main results (Theorems \ref{thm:connections-rigidity-discr} and
\ref{thm:connections-rigidity-cont}) we show that for separable Hilbert spaces the set of all
operators 
($C_0$-semigroups) having such 
behaviour is a residual set in an appropriate sense. 
This generalises
\cite{eisner:2008}, \cite{eisner/sereny:2008} and
\cite{eisner/sereny:2008-2} and extends a result of 
Choksi, Nadkarni \cite{choksi/nadkarni:1986} for unitary operators, see Nadkarni \cite[Chapter
7]{nadkarni:1998}. 
Connections to ergodic and measure theory are discussed. 

\section{Preliminaries}

We now make precise what we mean by rigidity. 
\begin{defi}
A bounded operator $T$ on a Banach space $X$ is called \emph{rigid} if
$$
 \text{strong-}\lim_{j\to\infty}T^{n_j}=I \quad \text{for some
 subsequence } 
 \{n_j\}_{l=j}^\infty\subset \N. 
$$
\end{defi}
\noindent Note that one can assume in the above definition 
$\lim_{j\to\infty}n_j=\infty$.
(Indeed, if $n_j$ does not converge to $\infty$,
then $T^{n_0}=I$ for some $n_0$ implying $T^{nn_0}=I$ for every $n\in\N$.)
So, rigidity describes a certain asymptotic behaviour of the powers $T^n$. 

\begin{remark}\label{remark:rigidity-then-unitary}
Rigid operators have no non-trivial weakly stable orbit. In
particular, by the Foia\c{s}--Sz.-Nagy decomposition, see
Sz.-Nagy, Foia\c{s} \cite{szokefalvi/foias:1960} and Foguel \cite{foguel:1963}, rigid contractions
on Hilbert spaces are necessarily unitary. 
\end{remark}

As trivial examples of rigid operators take $T:=\alpha I$ for $\alpha\in\Gamma$, $\Gamma$
the unit circle. Moreover, arbitrary (countable) 
combinations of such operators are rigid as well, as the following proposition shows.
\begin{prop}\label{prop:rigidity-discr-spectrum}
Let $X$ be a separable Banach space and let $T\in\LLL(X)$ be power
bounded with discrete spectrum, i.e., satisfying
$$
  H=\overline{\emph{\lin}}\{x\in X:\ Tx=\lambda x \text{ for some } \lambda\in \Gamma\}.
$$
Then $T$ is rigid. 
\end{prop}
\begin{proof}
Since $X$ is separable, the strong operator topology is metrisable on
bounded sets of $\LLL(H)$ (take for example the metric
$d(T,S):=\sum_{j=0}^\infty \|Tz_j-Sz_j\|/(2^j\|z_j\|)$ for a
dense sequence $\{z_j\}_{j=1}^\infty\subset X\setminus\{0\}$).
So it suffices to show that $I$ belongs to the strong closure of
$\{T^n\}_{n\in\N}$. For $\veps>0$, $m\in \N$ and $x_1,\ldots, x_m\in
X$, we have to find $n\in\N$ such that $\|T^nx_j-x_j\|<\veps$ for every
$j=1,\ldots, m$. 

Assume first that each $x_j$ is an eigenvector corresponding to some
unimodular eigenvalue $\lambda_j$, and hence $\|T^n x_j-x_j\|=|\lambda_j^n-1|
\|x\|$. Consider the compact group $\Gamma^m$ and the rotation $\vphi:\Gamma^m\to\Gamma^m$ given by $\vphi(z):=a z$
for $a:=(\lambda_1,\ldots, \lambda_m)$. By a classical recurrence theorem, see e.g. Furstenberg \cite[Theorem 1.2]{furstenberg-book}, 
there exists $n$ such that
$|\vphi^n(\Eins)-\Eins|<\veps_1:=\frac{\veps}{\max_{j=1,\ldots,m}\|x_j\|}$, i.e.,
$$
  |\lambda_j^n-1|<\veps_1 \quad
\text{ for every } j=1,\ldots,m,
$$ 
implying $\|T^nx_j-x_j\|<\veps$ for every
$j=1,\ldots, m$.

Assume now $0\neq x_j\in\lin\{x\in X:\ Tx=\lambda x \text{ for some }
\lambda\in \Gamma\}$ for all $j$. Then we have
$x_j=\sum_{k=1}^{K}c_{jk}y_k$ for $K\in\N$, eigenvectors $y_k\in X$, and
$c_{jk}\in\C$, $k=1,\ldots, K$, $j=1,\ldots, m$. Take
$\veps_2:=\frac{\veps}{K\max_{j,k}|c_{jk}|}$. 
 By the above, there exists $n\in\N$ such that
$\|T^ny_k-y_k\|<\veps_2$ for every $k=1,\ldots, K$ and therefore
$$
  \|Tx_j-x_j\|\leq \sum_{k=1}^{K}|c_{jk}|\|Ty_k-y_k\| <
\veps
$$       
for every $j=1,\ldots, m$. The standard density argument covers the
case of arbitrary $x_j\in X$.
\end{proof}

Analogously, one defines $\lambda$-rigid operators by replacing $I$ by $\lambda I$ in the above definition.

\begin{defi}
Let $X$ be a Banach space, $T\in\LLL(X)$ and  $\lambda\in\Gamma$.
We call $T$  \emph{$\lambda$-rigid} if there exists a sequence
$\{n_j\}_{j=1}^\infty\subset \N$ 
such that
$$
  \text{strong-} \lim_{j\to\infty}T^{n_j}=\lambda I.
$$  
\end{defi}
\noindent Again one can choose the sequence $\{n_j\}_{j=1}^\infty$ to  
converge to $\infty$.  
We finally call $T$ $\Gamma$-rigid if $T$ is $\lambda$-rigid for every $\lambda\in\Gamma$.

\begin{remark} Since every $\lambda$-rigid operator is $\lambda^n$-rigid for every $n\in\N$, we see that $\lambda$-rigidity implies rigidity. Moreover, $\lambda$-rigidity is equivalent to 
$\Gamma$-rigidity whenever $\lambda$ is irrational, i.e., $\lambda\notin e^{2\pi i\Q}$. (We used the fact
that for irrational $\lambda$ the set $\{\lambda^n\}_{n=1}^\infty$ is
dense in $\Gamma$ and that limit sets are always closed.) 
\end{remark}

The simplest examples are again operators of the form $\lambda I$,
$|\lambda|=1$. Indeed, $T=\lambda I$ is $\lambda$-rigid, and it is  
$\Gamma$-rigid if and only if $\lambda$ is irrational.  


\smallskip

Analogously, one can define rigidity for strongly continuous
semigroups. 
\begin{defi}
A $C_0$-semigroup $(T(t)_{t\geq 0})$ on a Banach space is called
\emph{$\lambda$-rigid} for $\lambda\in\Gamma$ if there exists a
sequence $\{t_j\}_{j=1}^\infty\subset \R_+$ with
$\lim_{j\to\infty} t_j=\infty$ such that 
$$
  \text{strong-}\lim_{j\to\infty}T(t_j)=\lambda I.
$$
Semigroups which are $1$-rigid are called \emph{rigid}, and semigroups
rigid for every $\lambda\in\Gamma$ are called \emph{$\Gamma$-rigid}.
\end{defi}
As in the discrete case, $\lambda$-rigidity for some $\lambda$ implies rigidity, and $\lambda$-rigidity for some irrational $\lambda$ is equivalent to $\Gamma$-rigidity. 

\begin{remark}\label{remark:rigidity-unitary-group}
Again, rigid semigroups have no non-zero weakly
stable orbit. This implies by the Foia\c{s}--Sz.-Nagy--Foguel
decomposition, see
Sz.-Nagy, Foia\c{s} \cite{szokefalvi/foias:1960} and Foguel \cite{foguel:1963} that every rigid semigroup on
a Hilbert space is automatically unitary.  
\end{remark}

The simplest examples of rigid $C_0$-semigroups are given by
$T(t)=e^{iat}I$ for some $a\in\R$. In this case, $T(\cdot)$ is
automatically $\Gamma$-rigid whenever $a\neq 0$. Moreover, one has the
following continuous analogue of Proposition \ref{prop:rigidity-discr-spectrum}.

\begin{prop}
Let $X$ be a separable Banach space and let $T(\cdot)$ be a 
bounded $C_0$-semigroup with discrete spectrum, i.e., satisfying
$$
  H=\overline{\emph{\lin}}\{x\in X:\ T(t)x=e^{ita}x \text{ for
    some } a\in \R \text{ and all } t\geq 0\}.
$$
Then $T$ is rigid. 
\end{prop}

So rigidity becomes non-trivial for operators ($C_0$-semigroups) having no point spectrum on the unit
circle. By a version of the classical Jacobs--Glicks\-berg--\-de Leeuw theorem (see
e.g. Krengel  \cite[pp. 108--110]{krengel:1985}) and for, i.e., power bounded operators 
on reflexive Banach spaces, 
the absence of point spectrum on $\Gamma$ is
equivalent to 
$$
  \text{weak-}\lim_{j\to\infty}T^{n_j}=0 \text{ for some subsequence }
  \{n_j\} \text{ with density }1,  
$$
where the density of a set $M\subset\N$ is defined by 
$$
  d(M):=\lim_{n\to\infty}\frac{\#(M\cap\{1,\ldots,n\})}{n}
$$ 
whenever this limit exists.
An analogous assertion holds for $C_0$-semigroups as well, see
e.g. \cite{EFNS}.
We call operators and semigroups with this property \emph{almost
  weakly stable}. (For a survey on almost weak stability in the
continuous case see Eisner, Farkas, Nagel, Ser\'eny \cite{EFNS}).

Restricting ourselves to Hilbert spaces, we will see that almost weakly stable and $\Gamma$-rigid operators and semigroups are the rule and not just an exception. 

 
%

\section{Discrete case: powers of operators}\label{section:rigidity-discrete}

%
%
%
%

We now introduce the spaces of operators we will work with.  

Take a separable infinite-dimensional
Hilbert space $H$ and denote by $\U$ the space of all unitary operators on
$H$ endowed with the strong$^*$ topology, i.e., the topology defined
by the seminorms 
$$
  p_x(T):=\sqrt{\|Tx\|^2+\|T^*x\|^2}
$$ 
(for details on this
topology see e.g. Takesaki \cite[p. 68]{takesaki}). Convergence in
this topology is strong convergence of operators and their adjoints.
Then $\U$ is a complete metric space with respect to the metric  
\begin{equation*}
 d(U,V):= \sum_{j=1}^\infty \frac{\|Ux_j -Vx_j\| + \|U^*x_j -V^*x_j\|}{2^j \|x_j\|}\quad  \text{for } U,V\in \U,
\end{equation*}
where $\{x_j\}_{j=1}^\infty$ is some dense subset of $H\setminus
\{0\}$. Similarly, the space $\I$ of
all isometric operators on $H$ will be endowed with the strong operator
topology and then is a complete metric space with respect to 
\begin{equation*}
d(T,S):= \sum_{j=1}^\infty \frac{\|Tx_j -Sx_j\|}{2^j \|x_j\|}\quad \text{for } T,S \in \mathcal{I}.
\end{equation*}
Finally, we denote by $\CC$ the space of all contractions on $H$
endowed with the weak operator topology which is a complete metric
space for the metric 
\begin{equation*}
d(T,S):= \sum_{j,k=1}^\infty \frac{|\la (T-S)x_j,x_k\ra|}{2^j
  \|x_j\| \|x_k\|}\quad \text{for } T,S \in \mathcal{\CC}.
\end{equation*}

The following is a basic step of our construction.


 
\begin{thm}\label{thm:rigidity-I}
Let $H$ be a separable infinite-dimensional Hilbert space. 
The set 
$$
 M:=\{T: \ 
\lim_{j\to \infty}T^{n_j}
 =I \text{ strongly for some }n_j\to\infty\}
$$
is residual for the weak
operator topology in the set $\CC$ of all contractions on $H$. This set is also residual for the strong operator
topology in the set $\I$ of all
isometries and in the set $\U$ of all unitary operators for the strong$^*$ operator topology.  
\end{thm}
\begin{proof}
We begin with the isometric case. 


Let $\{x_l\}_{l=1}^\infty$ be a dense subset of
$H\setminus\{0\}$. Since one can remove the assumption $\lim_{j\to\infty}n_j=\infty$ from the definition of
$M$, we have 
\begin{equation}\label{eq:rigidity-M-isom}
  M=\{T\in \I:\ \exists \{n_j\}_{j=1}^\infty\subset \N \text{ with } \lim_{j\to \infty} T^{n_j}x_l =x_l 
  \ \ \forall  l\in\N\}.
\end{equation}
Consider the sets 
$$
M_{k}:=\{T\in\I:\ \sum_{l=1}^\infty\frac{1}{2^l \|x_l\|}
\|T^n x_l - x_l \| < \frac{1}{k} \text{ for some }n\}
$$ 
which are open in the strong operator topology. 
Therefore, 
$$
 M=\bigcap_{k=1}^\infty M_k
$$ 
implies that $M$ is a $G_\delta$-set. To show that $M$ is residual  
it just remains to prove that $M$ is dense. Since $M$ contains all periodic unitary operators which are dense
in $\I$, see e.g. 
Eisner, Sereny \cite{eisner/sereny:2008}, the assertion follows. 

\smallskip

While for unitary operators the above arguments work as well, 
we need more delicate arguments for the space $\CC$ of all contractions. 

We first show that 
\begin{equation}\label{eq:rigidity-M}
  M=\{T\in \CC:\ \exists \{n_j\}\subset\N \text{ with } \lim_{j\to \infty}\la T^{n_j}x_l,x_l\ra =\|x_l\|^2 
  \ \ \forall  l\in\N\}.
\end{equation}
The inclusion ``$\subset$'' is clear. To prove the converse
inclusion, assume that $\lim_{j\to\infty}\la
T^{n_j}x_l,x_l\ra=\|x_l\|^2$ for each $l\in\N$. By the standard density argument we have
$\lim_{j\to\infty}\la T^{n_j}x,x\ra=\|x\|^2$ for every $x\in
H$. Strong convergence of $T^{n_j}$ to $I$ now follows from 
$$
   \|(T^{n_j}-I)x\|^2=\|T^{n_j}x\|^2 - 2\Re \la T^{n_j}x,x\ra + \|x\|^2\leq 2(\|x\|^2-\la T^{n_j}x,x\ra).
$$

We now define 
$$
M_{k}:=\{T\in\CC:\ \sum_{l=1}^\infty\frac{1}{2^l \|x_l\|^2}
|\la(T^n-I) x_l,x_l \ra| < \frac{1}{k} \text{ for some }n\},
$$ 
and observe again that $M=\bigcap_{k} M_{k}.$ 
 
It remains to show that the complement $M_{k}^c$ of $M_k$ is a nowhere dense
set. 
Since the set of periodic unitary operators $U_{per}$ on $H$ is
dense in the set of all contractions for the weak operator topology
(see e.g. Eisner, Ser\'eny \cite{eisner/sereny:2008}), it suffices to show $U_{per}\cap
\overline{M_{k}^c}=\emptyset$.
Assume that this is not the case,
i.e., that there exists a sequence $\{T_m\}_{m=1}^\infty\subset
M_{k}^c$ converging weakly to a periodic unitary operator $U$. 
Then by the standard argument (see, e.g., Eisner, Sereny \cite[Lemma 4.2]{eisner/sereny:2008}), $\lim_{m\to\infty}T_m=U$ strongly, hence $\lim_{m\to\infty}T_m^n=U^n$
strongly for every $n\in \N$. 
However, $T_m\in A_{k}^c$ means that
$$
  \sum_{l=1}^\infty \frac{1}{2^l \|x_l\|^2} |\la (T_m^n-I)x_l,x_l\ra| \geq
  \frac{1}{k} \quad \text{ for every }n,m\in\N.
$$ 
Since $T_m^n$ converges strongly and hence weakly to $U^n$ for every
$n$, and hence the expression on the left hand side of the above
inequality is dominated by $\{\frac{1}{2^{l-1}}\}$ which is a sequence in $l^1$, we
obtain by letting $m\to\infty$ that  
$$
  \sum_{l=1}^\infty \frac{1}{2^l \|x_l\|^2} |\la (U^n-I)x_l,x_l\ra| 
\geq
  \frac{1}{k} \quad \text{ for every }n,
$$ 
contradicting the periodicity of $U$.
\end{proof}

\smallskip

We now observe that one can replace $I$ in Theorem \ref{thm:rigidity-I} by $\lambda I$ for any $\lambda\in\Gamma$. To show this we need the following lemma. 

\begin{lemma}\label{lemma:rigidity-lambda-dense}
Let $H$ be a Hilbert space, $\lambda\in \Gamma$ and $N\in \N$. Then the set of all unitary operators $U$ with $U^n=\lambda I$ for some $n\geq N$ is dense in the set of all unitary operators for the norm topology.  
\end{lemma} 
\begin{proof}
Let $U$ be a unitary operator, $\lambda=e^{i\alpha}\in \Gamma$, $N\in
\N$ and $\veps>0$. By the spectral theorem 
$U$ is unitarily equivalent to a multiplication operator $\tilde{U}$
on some $L^2(\Omega, \mu)$ with 
\begin{equation*}
(\tilde{U}f)(\omega)=\varphi(\omega)f(\omega),\ \ \forall \omega\in\Omega,
\end{equation*} 
for some measurable $\varphi:\Omega \to \Gamma:=\{z\in\C:\ |z|=1\}$. 

We now approximate the operator $\tilde{U}$ as follows. Take $n\geq N$ such that $|1-e^\frac{2\pi i}{n}|\leq \veps$ and 
define for $\alpha_j:=e^{i\left( \frac{\alpha}{n}+ \frac{2\pi  j}{n}\right)}$, $j=0,\ldots, n-1$,
\begin{equation*}
\psi(\omega):= \alpha_{j-1},\ \forall \omega\in \varphi^{-1}(\{z\in\Gamma:\ \arg(\alpha_{j-1})\leq \arg(z)< \arg(\alpha_j) \}). 
\end{equation*} 
The multiplication operator $\tilde{P}$ corresponding to $\psi$ satisfies $\tilde{P}^{n}=e^{i\alpha}$. Moreover, 
\begin{equation*}
\Vert \tilde{U} - \tilde{P} \Vert = \sup_{\omega\in\Omega} |\varphi(\omega) - \psi(\omega)|\leq \varepsilon 
\end{equation*} 
proving the assertion.
\end{proof}

We now describe the ``typical'' asymptotic behaviour of  contractions
(iso\-metries, unitary operators) on separable Hilbert spaces. For an
alternative proof in the unitary case based on the spectral theorem
and an analogous result for measures on $\Gamma$ see Nadkarni \cite[Chapter 7]{nadkarni:1998}.  

\begin{thm}\label{thm:connections-rigidity-discr}
Let $H$ be a separable infinite-dimensional Hilbert space and $\Lambda\subset \Gamma$ be countable. Then the set of all operators $T$ satisfying the following properties
\begin{enumerate}
\item there exists $\{n_j\}_{j=1}^\infty\subset \N$ with density $1$ such that 
$$\lim_{j\to\infty}T^{n_j}=0 \quad \text{weakly},$$
\item for every $\lambda \in \Lambda$ there exists
  $\{n^{(\lambda)}_j\}_{j=1}^\infty$ with
  $\lim_{j\to\infty}n^{(\lambda)}_j=\infty$ such that 
$$\lim_{j\to\infty}T^{n^{(\lambda)}_j}=\lambda I \quad \text{strongly}$$
\end{enumerate}
is residual for the weak operator
topology in the set $\CC$ of all contractions. This set is also residual for the strong operator topology in the set $\I$ of all isometries as well as for the strong* operator topology in the set $\U$ of all unitary operators.  
\end{thm}
\noindent Recall that every contraction satisfying (2) is unitary,
cf. Remark \ref{remark:rigidity-then-unitary}.
%
%
\begin{proof}
By 
Eisner, Ser\'eny \cite{eisner/sereny:2008},
operators satisfying (1) 
are residual in $\CC$, $\I$ and $\U$. 

We now show that 
for a fixed $\lambda\in \Gamma$, the
set $M$ of all operators $T$ satisfying
$\lim_{j\to\infty}T^{n_j}=\lambda I$ strongly for some sequence
$\{n_j\}_{j=1}^\infty$ is residual. We again prove this first for
isometries and the strong operator topology. 

Take a dense set $\{x_l\}_{l=1}^\infty$ of $H\setminus \{0\}$ and observe that
$$
  M=\{T\in\I:\ \exists \{n_j\} \text{ with } \lim_{j\to
  \infty}T^{n_j}x_l =\lambda x_l 
  \ \forall  l\in\N\}.
$$
We see that $M=\bigcap_{k=1}^\infty M_k$ for the sets 
$$
  M_k:=\{T\in\I:\ \sum_{l=1}^\infty\frac{1}{2^l \|x_l\|}
\|T^n x_l-\lambda x_l\| < \frac{1}{k} \text{ for some }n\}
$$
which are open for the strong operator topology. Therefore $M$ is a
$G_\delta$-set which is dense by Lemma
\ref{lemma:rigidity-lambda-dense}, and the residuality of $M$ follows.

\smallskip

The unitary case goes analogously, and we now prove the more delicate
contraction case. To do so we first show that 
$$
  M=\{T\in\CC:\ \exists \{n_j\}_{j=1}^\infty \text{ with } \lim_{j\to \infty}\la T^{n_j}x_l,x_l\ra =\lambda\|x_l\|^2 
  \ \forall  l\in\N\}.
$$
As in the proof of Theorem \ref{thm:rigidity-I}, the nontrivial inclusion follows from 
\begin{eqnarray*}
  \|(T^{n_j}-\lambda I)x\|^2&=&\|T^{n_j}x\|^2 - 2\Re \la T^{n_j}x, \lambda x\ra + \|x\|^2\\
&\leq& 2(\|x\|^2-\la T^{n_j}x,\lambda x\ra)= 2\overline{\lambda}\la(\lambda I -T^{n_j})x, x\ra.
\end{eqnarray*}
For the sets
$$
  M_{k}:=\{T\in\CC:\ \sum_{l=1}^\infty\frac{1}{2^l}
|\la(T^n-\lambda I) x_l,x_l \ra| < \frac{1}{k} \text{ for some }n\},
$$
we have the equality $M=\bigcap_{k=1}^\infty M_k$. 
Note again that it is not clear whether
the sets $M_k$ are open for the weak operator topology, so we use another
argument to show that the complements $M_k^c$ are nowhere dense.  
By Lemma \ref{lemma:rigidity-lambda-dense} it suffices to show that
$M_k^c \cap U_\lambda=\emptyset$ for the complement $M_k^c$ and the
set $U_\lambda$ of all unitary operators $U$ satisfying $U^n=\lambda
I$ for some $n\in\N$. This can be shown analogously to the proof of Theorem \ref{thm:rigidity-I} by replacing $I$ by $\lambda I$. 


\smallskip

Since $\lambda$-rigidity for an irrational $\lambda$ already implies
$\Gamma$-rigidity, the proof is complete.
%
%
\end{proof}
%
%

We now present basic constructions leading to examples of operators
with properties described in Theorem \ref{thm:connections-rigidity-discr}.

\begin{example}\label{ex:rigidity-from-measures}
a) A large class of abstract examples of $\Gamma$-rigid unitary operators which are
almost weakly stable comes from harmonic analysis. There, a
probability measure $\mu$ on $\Gamma$ is called $\lambda$-rigid if its
Fourier coefficients satisfy
$$
  \lim_{j\to\infty}\hat{\mu}_{n_j}=\lambda\quad \text{for some }
  \{n_j\}_{j=1}^\infty\subset\N, \ \lim n_j=\infty. 
$$
A result of Choksi, Nadkarni \cite{choksi/nadkarni:1986}, see Nadkarni \cite[Chapter
  7]{nadkarni:1998} states that $\lambda$-rigid continuous measures
form a dense $G_\delta$ (and hence a residual) set in the space of all
probability measures with respect to the weak$^*$ topology. 

Take a $\lambda$-rigid measure $\mu$ for some irrational $\lambda$. 
Note that the arguments used for the existence of such a measure are based on the Baire category theorem
without yielding any concrete example. 
The unitary operator given by $(Uf)(z):=zf(z)$ on $L^2(\Gamma,\mu)$ satisfies conditions (1) and (2) of Theorem \ref{thm:connections-rigidity-discr}. 
To show (2), it again suffices to prove that $U$ is $\lambda$-rigid
since $\lambda$ is irrational. By our assumption on $\mu$, there exists a subsequence $\{n_k\}_{k=1}^\infty$ such that the Fourier
coefficients of $\mu$ satisfy
$\lim_{k\to\infty}\hat{\mu}_{n_k}=\lambda$. A standard argument based
on the fact that the characters $\{z\mapsto z^n\}_{n=-\infty}^\infty$ form an orthonormal basis in
$L^2(\Gamma,\mu)$ implies that $\lim_{k\to\infty}U^{n_k}=\lambda I$ weakly, and
hence strongly. 


Conversely, if $U$ satisfies conditions (1) and (2) of Theorem \ref{thm:connections-rigidity-discr},
then every spectral measure $\mu$ of $U$ is continuous and
$\lambda$-rigid for every $\lambda\in\Gamma$.

\smallskip

b) Another class of examples of rigid 
almost weakly stable unitary operators
comes from ergodic theory. 
For a probability space $(\Omega, \mu)$ and a measurable,
$\mu$-preserving transformation $\vphi:\Omega\to \Omega$, the
associated unitary operator $U$ on $H:=L^2(\Omega, \mu)$ is given by
$(Uf)(\omega):=f(\vphi(\omega))$. 
The transformation $\vphi$ is called rigid if $U$ is rigid, and
$\lambda$-rigid or $\lambda$-weakly mixing if the restriction $U_0$ of $U$ to the invariant subspace
$H_0:=\{f: \int_\Omega f d\mu =0\}$ is $\lambda$-rigid. 
This restriction $U_0$ is almost weakly stable, i.e.,
satisfies (2) in Theorem \ref{thm:connections-rigidity-discr} if and
only if $\vphi$ is weakly mixing. Thus, each weakly mixing rigid (or
$\lambda$-rigid) transformation corresponds to a unitary operator
satisfying conditions (1) and (2) in Theorem
\ref{thm:connections-rigidity-discr}.

Katok \cite{katok:1985} proved that rigid transformations form a dense
$G_\delta$-set in the set of all measure preserving transformations,
and Choksi, Nadkarni \cite{choksi/nadkarni:1986} generalised this result to $\lambda$-rigid transformations.  
For more information we refer to Nadkarni
\cite[p. 59]{nadkarni:1998} and for concrete examples of rigid weakly mixing transformations using adding machines and interval exchange transformations see Goodson, Kwiatkowski, Lemanczyk, Liardet \cite{lemanczyk-etal:1992} and Ferenczi, Holton, Zamboni \cite{ferenczi/holton/zamboni:2005}, respectively. For examples of rigid weakly mixing transformations given by
Gaussian automorphisms see Cornfeld, Fomin, Sinai \cite[Chapter 14]{cornfeld/fomin/sinai:1982}.  
\end{example}

Furthermore, there is an (abstract) method of constructing
$\lambda$-rigid operators from a rigid one. The idea of this
construction in the context of measures belongs to Nadkarni
\cite[Chapter 7]{nadkarni:1998}.

\begin{example}\label{ex:rigidity-rescaling}
Let $T$ be a rigid contraction with $\lim_{j\to\infty}T^{n_j}=I$
strongly, and let $\lambda\in \Gamma$. We construct a class of
$\lambda$-rigid operators from $T$. Note that if $\lambda$ is
irrational and if $T$ is unitary with no point
spectrum, this construction gives us a class of examples satisfying (1)-(2) of Theorem \ref{thm:connections-rigidity-discr}.  

Take $\alpha\in\Gamma$ and consider the operator $T_\alpha:=\alpha
T$. Then we see that $T_{\alpha}$ is $\lambda$-rigid if
$\lim_{j\to\infty}\alpha^{n_j}=\lambda$ for the above sequence
$\{n_j\}$. Nadkarni \cite[p. 49--50]{nadkarni:1998} showed that the set of all $\alpha$ such
that the limit set of $\{\alpha^{n_j}\}_{j=1}^\infty$ contains an
irrational number has full Lebesgue measure in
$\Gamma$. Every such $\alpha$ leads to a $\Gamma$-rigid
operator $\alpha T$.  
%
\end{example}
%

We now show that one cannot replace the operators $\lambda I$
in Theorem \ref{thm:connections-rigidity-discr} by any other operator.

\begin{prop}\label{prop:rigidity-limit-is-const}
Let $V\in\LLL(H)$ be such that the set 
$$
   M_V:=\{T: \exists \{n_j\}_{j=1}^\infty \text{ such that } \lim_{j\to \infty}T^{n_j} =V \text{ strongly}\}
$$
is dense in one of the spaces $\U$, $\I$ or $\CC$. Then $V$ is a
multiple of identity.
\end{prop}
\begin{proof}
Consider the contraction case and assume that the set of all
contractions $T$ such that weak-$\lim_{j\to\infty}T^{n_j}=V$ for some $\{n_j\}_{j=1}^\infty$ is dense in $\CC$. Since
every such operator $T$ commutes with $V$ by
$TV=\lim_{j\to\infty}T^{n_j+1}=VT$, we obtain by assumption that $V$
commutes with every contraction. In particular, $V$ commutes with
every one-dimensional projection implying that 
$V=\lambda I$ for some $\lambda\in\C$.  

The same argument works for the spaces $\I$ and $\U$ using the
density of unitary operators in the set of all contractions for the
weak operator topology.
\end{proof}
\begin{remark}
In the above proposition, one has $V=\lambda I$ for some $\lambda\in \Gamma$
in the unitary and isometric case. Moreover, the same holds in the
contraction case if $M_V$ is residual. 
(This follows from 
the fact that 
the set of all non-unitary 
contractions is of first category in $\CC$ by Remark
\ref{remark:rigidity-then-unitary} and Theorem \ref{thm:rigidity-I},
see also \cite{eisner:2008}).
\end{remark}

\begin{remark}
It is not clear whether Theorem \ref{thm:connections-rigidity-discr}
remains valid under the additional requirement 
\begin{equation}\label{eq:rigidity-inside-Gamma}
\{\lambda\in\C:|\lambda|<1\}\cdot I\subset \ol{\{T^n:n\in\N\}}^\sigma, 
\end{equation}
where $\sigma$ denotes the weak operator topology. 
Since countable intersections of residual sets are residual and 
the right hand side of (\ref{eq:rigidity-inside-Gamma}) is closed, this question becomes whether, for a fixed $\lambda$ with $0<|\lambda|<1$, the set
$M_\lambda$ of all contractions $T$ satisfying $\lambda I\in
\ol{\{T^n:n\in\N\}}^\sigma$ is residual.
Note that each $M_\lambda$ is dense in $\CC$ since, for $\lambda=re^{is}$, it contains the
set 
$$
  \{cU: 0<c<1,\ c^n=r\text{ and }U^n=e^{is}I \text{ for some }n\}
$$
which is dense in $\CC$ for the
norm topology by Lemma \ref{lemma:rigidity-lambda-dense}.
%
\end{remark}

We finally mention that absence of rigidity does not imply weak
stability, or, equivalently, absence of weak stability does not imply rigidity, as
the following example shows. 


\begin{example}
There exist 
unitary operators $T$ with no non-trivial weakly stable orbit which are nowhere rigid, i.e., such
that $\lim_{j\to\infty}T^{n_j}x=x$ for some subsequence
$\{n_j\}_{j=1}^\infty$ implies $x=0$. (Note that such operators are automatically 
almost weakly stable by Proposition \ref{prop:rigidity-discr-spectrum}.) A class of such examples
comes from \emph{mildly mixing} transformations which are not strongly mixing, see Furstenberg, Weiss \cite{furstenberg/weiss:1978}, Fra\c{c}zek, Lema\'nczyk \cite{fraczek/lemanczyk:2006} and Fra\c{c}zek, Lema\'nczyk, Lesigne \cite{fraczek/lemanczyk/lesigne:2007}.
\end{example}


\section{Continuous case: $C_0$-semigroups} 

We now give the continuous analogue of the above results for unitary and isometric strongly continuous
(semi)groups. 
 
Let $H$ be again a separable infinite-dimensional Hilbert space. We
denote by $\U^{\cont}$ the set of all unitary $C_0$-groups on $H$
endowed with the topology of strong convergence of semigroups and their adjoints uniformly on compact time intervals.
This is a complete metric and hence a Baire space for  
\begin{equation*}
d(U(\cdot),V(\cdot)):= \sum_{n,j=1}^\infty
\frac{\sup_{t\in[-n,n]}\|U(t)x_j -V(t)x_j\|}
{2^j \|x_j\|}\quad  \text{for } U(\cdot),V(\cdot)\in \U^{\cont},
\end{equation*}
where $\{x_j\}_{j=1}^\infty$ is a fixed dense subset of $H\setminus\{0\}$.
We further denote by $\I^{\cont}$ the set of all isometric $C_0$-semigroups on $H$
endowed with the topology of strong convergence uniform on compact
time intervals. Again, this is a complete metric space for  
 \begin{equation*}
d(T(\cdot),S(\cdot)):= \sum_{n,j=1}^\infty \frac{\sup_{t\in[0,n]}\|T(t)x_j -S(t)x_j\|}{2^j \|x_j\|}\quad \text{for } T(\cdot),S(\cdot) \in \mathcal{I}^{\cont}.
\end{equation*}


The proofs of the following results are similar to the discrete case, but require some additional technical details.


\begin{thm}\label{thm:rigidity-I-cont}
Let $H$ be a separable infinite-dimensional Hilbert space. 
The set 
$$
 M^{cont}:=\{T(\cdot): \  \lim_{j\to \infty}T(t_j)=I \text{ strongly
   for some } t_j\to\infty\}
$$
is residual in the set $\I^{\cont}$ of all isometric $C_0$-semigroups for the
topology corresponding to strong convergence uniform on compact time
intervals in $\R_+$. The same holds for unitary $C_0$-groups and the topology of strong convergence uniform on compact time intervals in $\R$. 
\end{thm}
\begin{proof}
We begin with the unitary case. 

Choose $\{x_l\}_{l=1}^\infty$ as a dense subset of
$H\setminus\{0\}$. Since one can replace $\lim_{j\to\infty}t_j=\infty$
in the definition of $M^{\cont}$ by $\{t_j\}_{j=1}^\infty\subset
[1,\infty)$ we have 
\begin{equation}\label{eq:rigidity-M-cont}
  M^{\cont}=\{T(\cdot)\in\U^{\cont}:\ \exists \{t_j\}\in [1,\infty): \ \  
\lim_{j\to \infty}T(t_j)x_l =x_l 
  \ \ \forall  l\in\N\}.
\end{equation}
%
%
Consider now the open sets 
$$
M_{k,t}:=\{T(\cdot)\in\U^{\cont}:\ \sum_{l=1}^\infty\frac{1}{2^l \|x_l\|}
\|T(t) x_l- x_l \| < \frac{1}{k}\} 
$$ 
and $M_k^{\cont}:=\bigcup_{t\geq 1} M_{k,t}$. We have 
$$
  M^{\cont}=\bigcap_{k=1}^\infty M_{k}^{\cont},
$$
and hence $M^{\cont}$ is a $G_\delta$-set.
%
Since periodic unitary $C_0$-groups are dense in  $\U^{\cont}$ by
Eisner, Ser\'eny \cite{eisner/sereny:2009}, and since they are
contained in $M^{\cont}$, we see that $M^{\cont}$ is residual as a countable intersection of dense open sets. 

The same arguments and the density of periodic unitary operators in
$\I^{\cont}$ (see Eisner, Ser\'eny \cite{eisner/sereny:2009}) imply the assertion in the isometric case.
\end{proof}

The following continuous analogue of Lemma
\ref{lemma:rigidity-lambda-dense} allows to replace $I$ by $\lambda I$. 
\begin{lemma}\label{lemma:rigidity-lambda-dense-cont}
Let $H$ be a Hilbert space and fix $\lambda\in \Gamma$ and $N\in
\N$. Then for every unitary $C_0$-group $U(\cdot)$ there exists a
sequence $\{U_n(\cdot)\}_{n=1}^\infty$ of unitary $C_0$-groups such
that 
\begin{enumerate}[(a)]
\item For every $n\in\N$ there exists $\tau\geq N$ with $U_n(\tau)=\lambda I$, 
\item $\lim_{n\to\infty}\|U_n(t)-U(t)\|=0$ uniformly on compact
  intervals in $\R$.
\end{enumerate}
\end{lemma} 
\begin{proof}
Let $U(\cdot)$ be a unitary $C_0$-group on $H$,
$\lambda=e^{i\alpha}\in \Gamma$. By the
spectral theorem, $H$ is isomorphic to $L^2(\Omega, \mu)$ for some
finite measure space $(\Omega, \mu)$ and $U(\cdot)$ is unitarily
equivalent to a multiplication group $\tilde{U}(\cdot)$ given by 
\begin{equation*}
(\tilde{U}(t)f)(\omega)=e^{itq(\omega)}f(\omega),\quad  \omega\in\Omega,
\end{equation*} 
for some measurable $q:\Omega \to \R$. 

To approximate $\tilde{U}(\cdot)$, let $N\in \N$, $\veps>0$, $t_0>0$ and take
$m\geq N$, $m\in \N$, such that $\|1-e^\frac{2\pi i}{m}\|\leq \veps/(2 t_0)$. Define for $\alpha_j:=e^{i\left( \frac{\alpha}{m}+ \frac{2\pi  j}{m}\right)}$, $j=0,\ldots, m-1$,
\begin{equation*}
p(\omega):= \alpha_{j-1} \ \text{ for all } \omega\in \varphi^{-1}(\{z\in\Gamma:\ \arg(\alpha_{j-1})\leq \arg(z)< \arg(\alpha_j) \}). 
\end{equation*} 
The multiplication group $\tilde{V}(\cdot)$ defined by
$\tilde{V}(t)f(\omega):=e^{itp(\omega)f(\omega)}$ satisfies
$\tilde{V}(m)=e^{i\alpha}$. Moreover, 
%
\begin{eqnarray*}
\Vert \tilde{U}(t)f - \tilde{V}(t)f \Vert^2 &=& \int_\Omega
|e^{itq(\omega)} - e^{itp(\omega)}|^2 \|f(\omega)\|^2
\\ 
&\leq& 2|t|\sup_{\omega\in\Omega}|q(\omega)-p(\omega)| \|f\|^2 <
\veps \|f\|^2  
\end{eqnarray*}
uniformly in $t\in[-t_0,t_0]$. 
\end{proof}
%

We now
obtain the following characterisation of the ``typical'' asymptotic
behaviour of isometric and unitary $C_0$-(semi)groups on separable
Hilbert spaces.

\begin{thm}\label{thm:connections-rigidity-cont}
Let $H$ be a separable infinite-dimensional Hilbert space. Then the
set of all $C_0$-semigroups $T(\cdot)$ on $H$ satisfying the following properties
\begin{enumerate}
\item there exists a set $M\subset \R_+$ with density $1$ such that 
$$
  \lim_{t\to\infty, t\in M}T(t)=0 \quad \text{weakly},
$$
\item for every $\lambda \in \Gamma$ there exists
  $\{t^{(\lambda)}_j\}_{j=1}^\infty$ with
  $\lim_{j\to\infty}t^{(\lambda)}_j=\infty$ such that 
$$
  \lim_{j\to\infty}T(t^{(\lambda)}_j)=\lambda I \quad \text{strongly}
$$
\end{enumerate}
is residual in the set of all isometric $C_0$-semigroups for the
topology of strong convergence uniform on compact time intervals in $\R_+$. The
same holds for unitary $C_0$-groups for the topology of strong convergence uniform on compact time intervals in $\R$.  
\end{thm}
\noindent Recall that the density of a
measurable set $M\subset \R_+$ is 
$$
  d(M):=\lim_{t\to\infty}\frac{\mu(M\cap [0,t])}{t}\leq 1 
$$
whenever the limit exists. 
\begin{proof}
By  
Eisner, Ser\'eny \cite{eisner/sereny:2009},
$C_0$-(semi)groups satisfying (1)
are residual in $\I^\cont$ and $\U^\cont$. 

We show, for a fixed $\lambda\in \Gamma$, the
residuality of the set $M^{(\lambda)}$ of all $C_0$-semigroups $T(\cdot)$ satisfying
strong-$\lim_{j\to\infty}T(t_j)=\lambda I$ for some sequence
$\{t_j\}_{j=1}^\infty$ converging to infinity. We prove this property for
the space $\I^{\cont}$ of all isometric semigroups and the strong
operator convergence uniform on compact intervals, the unitary case
goes analogously. 

Take $\lambda\in\Gamma$ and observe 
$$ 
  M^{(\lambda)}=\{T(\cdot)\in\I^{\cont}:\ \exists t_j\to\infty \text{
  with } \lim_{j\to \infty} T(t_j)x_l =\lambda x_l \ \forall  l\in\N\}
$$
for a fixed dense sequence $\{x_l\}_{l=1}^\infty\subset H\setminus\{0\}$.
%
Consider now the open sets  
$$
  M_{k,t}:=\left\{T(\cdot)\in\I^{\cont}:\ \sum_{l=1}^\infty\frac{\|(T(t)-\lambda I) x_l \|}{2^l \|x_l\|} < \frac{1}{k}\right\} 
$$
and their union $M_k:=\bigcup_{t\geq 1}M_{k,t}$ being open
as well. The equality $M^{(\lambda)}=\bigcap_{k=1}^\infty M_k$ follows
as in the proof of Theorem \ref{thm:rigidity-I-cont}. Since every $M_k$ contains
periodic unitary $C_0$-groups and is therefore dense by Lemma
\ref{lemma:rigidity-lambda-dense-cont}, $M^{(\lambda)}$ is residual as
a dense countable intersection of open sets. 

Since 
$\lambda$-rigidity for some irrational $\lambda$ already implies $\lambda$-rigidity for every $\lambda\in\Gamma$,
the theorem is proved.
%
%
\end{proof}
Note that every semigroup satisfying (1) and (2) above is a unitary group by Remark \ref{remark:rigidity-unitary-group}.
\begin{example} 
a) There is the same correspondence between rigid (or $\lambda$-rigid)
unitary $C_0$-groups and rigid (or $\lambda$-rigid) probability
measures on $\R$ as in the discrete case, see Example
\ref{ex:rigidity-from-measures}. Here, to a probability measure $\mu$
on $\R$ one associates the multiplication group given by $(T(t))f(s):=e^{ist}f(s)$ on $H=L^2(\R,\mu)$. For $\lambda\in\Gamma$, we call a measure $\mu$ on $\R$ \emph{$\lambda$-rigid} if there exists $t_j\to\infty$, $t_j\in\R$, such
that the Fourier transform of $\mu$ satisfies
$\lim_{j\to\infty}\mathcal{F}\mu(t_j)=\lambda$. Using exactly the same
arguments as in Choksi, Nadkarni \cite{choksi/nadkarni:1986}, see
Nadkarni \cite[Chapter 7]{nadkarni:1998}, one shows that the set of
all continuous $\Gamma$-rigid measures on $\R$ is a dense $G_\delta$
set in the set of all Radon measures with respect to the weak$^*$
topology. For each such measure, the associated unitary $C_0$-group is
$\Gamma$-rigid and almost weakly stable, and conversely, the spectral
measures of an almost weakly stable $\Gamma$-rigid unitary group is
continuous and $\Gamma$-rigid.  
 

b) Again, another large class of examples comes from
ergodic theory. Consider a \emph{measure preserving semiflow}, i.e., a family of measure preserving
transformations $\{\vphi_t\}_{t\geq \R_+}$ on a probability space
$(\Omega, \mu)$ such that the function $(t,\omega)\mapsto
\vphi_t(\omega)$ is measurable on $\R_+\times\Omega$. 
On the Hilbert space $H:=L^2(\Omega, \mu)$ the semiflow induces an isometric semigroup
 by $(T(t))f(\omega):=f(\vphi_t(\omega))$ which is strongly continuous by Krengel \cite[{\S}1.6, Thm. 6.13]{krengel:1985}. 
If one/every $\vphi_t$ is invertible, i.e., if we start by a \emph{flow},
then $T(\cdot)$ extends to a unitary group.  

The semigroup $T(\cdot)$ is almost weakly stable and
$\lambda$-rigid if and only if the semiflow $(\vphi_t)$ is weakly
mixing and $\lambda$-rigid, where the last notion is defined
analogously to the discrete case. However, such flows are not so
well-studied and there seems to be no concrete construction for flows
in ergodic theory as it was done for operators. So we apply the
following abstract argument to present a large class of such flows. 

We start from a 
probability space $(\Omega, \mu)$. As discussed above, a ``typical''
measure preserving transformation $\vphi$ on $(\Omega,\mu)$ is weakly
mixing and $\Gamma$-rigid. On the other hand, by de la Rue, de Sam
Lazaro  \cite{delarue/desamlazaro:2003} a ``typical'' (for the same topology) measure
preserving transformation $\vphi$ is embeddable into a flow, i.e.,
there exists a flow $(\vphi_t)_{t\in\R}$ such that
$\vphi=\vphi_1$. Therefore, a ``typical'' transformation is weakly
mixing, $\Gamma$-rigid \emph{and} embeddable, and every such
transformation leads to an almost weakly stable $\Gamma$-rigid unitary group.  

c) Analogously to b), we can construct a class of examples 
on arbitrary Hilbert spaces. We use that every unitary operator $T$ is embeddable
into a unitary $C_0$-group $T(\cdot)$, i.e., $T(1)=T$ for some $T(\cdot)$, see e.g. \cite{eisner:2008-2}. Take now any operator
satisfying assertions of Theorem
\ref{thm:connections-rigidity-discr}. Since such an operator is
automatically unitary by Remark \ref{remark:rigidity-then-unitary}, it is
embeddable. Thus, every such operator leads to an example of a
$C_0$-group satisfying (1) and (2) of Theorem
\ref{thm:connections-rigidity-cont}. (Note that condition (1) follows
from the spectral mapping theorem for the point spectrum, see
e.g. Engel, Nagel \cite[Theorem IV.3.7]{engel/nagel:2000}.)
%
%
%
\end{example}
We again show that limit operators $\lambda I$, $|\lambda|=1$, cannot
be replaced in the above theorem by any other operator.   

\begin{prop}
Let for some $V\in\LLL(H)$ the set 
$$
   M^{\cont}_V:=\{T(\cdot): \exists t_j\to\infty \text{ such that } \lim_{j\to \infty}T(t_j)=V \text{ strongly}\}
$$ 
be dense in $\I^{\cont}$ or $\U^{\cont}$. Then $V=\lambda I$ for some $\lambda\in\Gamma$.
\end{prop}
\begin{proof} We prove this assertion for $\I$, the unitary case is
  analogous.

Observe that $V$ commutes with every $T(\cdot)\in M^{\cont}_V$ by 
$$
 VT(t)=\text{strong-}\lim_{j\to\infty} T(t_j+t)=T(t)V
$$ 
implying that $V$ commutes with every unitary 
$C_0$-group. Since unitary $C_0$-groups are dense in the set of all
contractive $C_0$-semigroups for the topology of weak operator convergence
uniform on compact time intervals by Krol \cite{krol:2009}, we see that $V$ commutes with every contractive
$C_0$-semigroup. We observe now that orthogonal one-dimensional projections are
embeddable into a contractive $C_0$-semigroup by
\cite[Prop. 4.7]{eisner:2008-2} and its proof. So $V$ commutes with
every orthogonal one-dimensional projection, which implies $V=\lambda I$ for some $\lambda\in\C$. 
Moreover, $|\lambda|=1$ holds since $V$ is the strong limit of isometric operators. 
\end{proof}

\begin{remark}
Recall that the space of all contractive $C_0$-semigroups on $H$ is
neither complete metric nor compact with respect to the  
topology of weak operator convergence uniform on compact time
intervals, see Eisner, Ser\'eny \cite{eisner/sereny:2008-2}.
So it is not clear whether one can
formulate an analogue of the above result for contractive
$C_0$-semigroups as done in the discrete case.
\end{remark}


\section{Further remarks}

We now consider some generalisations of the above results. 

\smallskip

\noindent \textbf{``Controlling'' the sequences $\{n_j\}$ and $\{t_j\}$.} 
We now take a closer look at the sequences $\{n_j\}$ and $\{t_j\}$ occuring in Theorems \ref{thm:connections-rigidity-discr}(2) and \ref{thm:connections-rigidity-cont}(2).   

\smallskip

Observe first that, by the same arguments as in the proofs of Theorems \ref{thm:rigidity-I}
and \ref{thm:connections-rigidity-discr}, we can replace $T$ there by $T^{m}$ for a fixed $m$. Changing
appropriately the assertion and the proof of Lemma \ref{lemma:rigidity-lambda-dense}, we see that 
one can add the condition $\{n_j\}_{j=1}^\infty\subset m\N$ to the
sequence appearing in rigidity and $\lambda$-rigidity. More precisely, for every
$\lambda\in \Gamma$ and $m\in \N$, the set of all operators $T$ such
that strong-$\lim_{j\to\infty}T^{n_j}=I$ for some $\{n_j\}\subset
m\N$ is residual in $\U$, $\I$ and $\CC$. 

\smallskip 

It is a hard problem to determine the sequences $\{n_j\}$ and $\{t_j\}$ exactly.  
However, one can generalise the above observation and ``control'' these sequences in the following sense. Let $\Lambda\subset \N$ be an unbounded set. We call an operator $T$ \emph{rigid along $\Lambda$} if strong-$\lim_{j\to\infty} T^{n_j}=I$ for some increasing sequence $\{n_j\}_{j=1}^\infty\subset \Lambda$. Similarly, we define rigidity along an unbounded set $\Lambda\subset\R_+$ for $C_0$-semigroups, $\lambda$- and $\Gamma$-rigidity. 
It follows from a natural modification of Lemmas \ref{lemma:rigidity-lambda-dense} and \ref{lemma:rigidity-lambda-dense-cont} 
that, for a fixed unbounded set $\Lambda$ in $\N$ and $\R_+$, respectively,  
one can assume $\{n^{(\lambda)}_j\}\subset \Lambda$ and $\{t^{(\lambda)}_j\}\subset \Lambda$ in Theorems \ref{thm:connections-rigidity-discr}(2) and \ref{thm:connections-rigidity-cont}(2). Thus, the set of all $\Gamma$-rigid operators (semigroups)  along a fixed unbounded set is residual in $\U$, $\I$ and $\CC$ ($\U^\cont$ and $\I^\cont$, respectively).


\smallskip

\noindent \textbf{Banach space case.}
We finally discuss briefly the situation in Banach spaces.

Note first that Theorems \ref{thm:connections-rigidity-discr} and
\ref{thm:connections-rigidity-cont} are not true in general separable
Banach spaces. Indeed, since weak convergence in $l^1$ implies strong
convergence, we see that (2) implies strong convergence to zero of
$T^n$ (or of $T(t)$, respectively), making (3) or
just rigidity impossible.   

We now consider the question in which Banach spaces 
rigid and $\Gamma$-rigid operators are residual. Since 
in the contraction case our techniques heavily use Hilbert space methods, we only consider the isometric and unitary case. 
Let $X$ be a separable infinite-dimensional Banach
space, and $\I$ be the set of all isometries on $X$ endowed with the
strong operator topology. Observe that the sets 
$$
M_{k}:=\{T\in\I:\ \sum_{l=1}^\infty\frac{1}{2^l \|x_l\|}
\|(T^n-I) x_l\| < \frac{1}{k} \text{ for some }n\}
$$  
appearing in the proof of Theorem \ref{thm:rigidity-I} for the
isometric case are still open, and therefore $M$ is a $G_\delta$-set containing
periodic isometries. Thus Theorem \ref{thm:rigidity-I} holds in all
separable infinite-dimensional Banach spaces such that periodic
isometries form a dense set of $\I$. Analogously, the set of operators satisying property (3) in
Theorem \ref{thm:connections-rigidity-discr} is residual in $\I$ if
and only if it is dense in $\I$. This is the case whenever there
exists an irrational $\lambda\in\Gamma$ such that the set of all isometries $T$ with
strong-$\lim_{j\to\infty}T^{n_j}=\lambda$ for some sequence $\{n_j\}$
is dense in $\I$. Analogously, the set of operators satisfying (1) and
(2) of Theorem \ref{thm:connections-rigidity-discr} is residual in $\I$ if
and only if it is dense in $\I$. The same assertions hold for the
set $\U$ of all invertible isometric operators with the topology induced by the seminorms $p_x(T)=\sqrt{\|Tx\|^2+\|T^{-1}x\|^2}$, which is a complete
metric space with respect to the metric
$$
  d(T,S)=\sum_{j=1}^\infty\frac{\|Tx_j-Sx_j\|+\|T^{-1}x_j-S^{-1}x_j\|}{2^j \|x_j\|}
$$ 
for a fixed dense sequence $\{x_j\}_{j=1}^\infty\subset X\setminus\{0\}$.

\bigskip

\smallskip

\noindent {\bf Acknowledgement.}
The author is grateful to Mariusz Lema\'nczyk, Mi\-cha\-el Lin, Vladimir
M\"uller, Rainer Nagel and Marco Schreiber for helpful comments, as well as Ralph Chill
and the University of Metz for the hospitality during her research visit. The work on this paper
was supported by the European Social Fund in Baden-W\"urttemberg.

\parindent0pt

\end{document}